\documentclass[a4paper,german,final,twoside,notitlepage,12pt]{article}
\usepackage[T1]{fontenc}
\usepackage{amssymb}
\usepackage{amsthm}
\usepackage{amsmath}
\usepackage[all]{xy}
\usepackage{graphicx}

\newtheorem{theorem}{Theorem}[section]
\newtheorem{definition}[theorem]{Definition}

\newtheorem{proposition}[theorem]{Proposition}
\newtheorem{lemma}[theorem]{Lemma}

\def\id {\mathrm{id}}

\newcommand{\alxy}[1]{\begin{aligned}\xymatrix{#1}\end{aligned}}
\newcommand{\alxydim}[2]{\begin{aligned}\xymatrix#1{#2}\end{aligned}}
\def\proof{\textit{Proof. }}
\def\endofproof{\hfill{$\square$}\\}

\numberwithin{theorem}{section}

\makeatletter

\def\adress#1{\gdef\@adress{#1}}
\def\@adress{}
\def\preprint#1{\gdef\@preprint{#1}}
\def\@preprint{}
\def\@maketitle{
  \newpage
  \noindent
  \begin{tabular}{cc}
    \begin{minipage}[c]{0.4\textwidth}
      \begin{flushleft}
        \includegraphics[width=110pt]{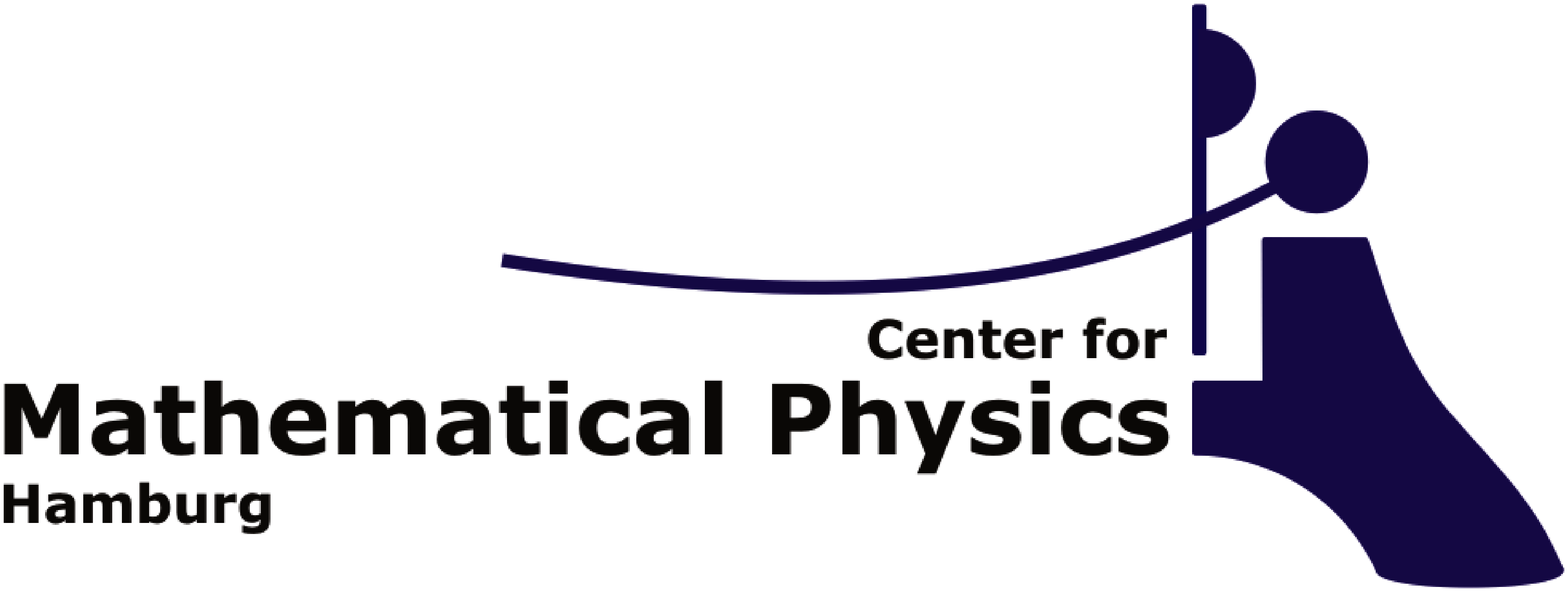}
      \end{flushleft}  
    \end{minipage}&
    \begin{minipage}[c]{0.55\textwidth}
      \begin{flushright}
      {\small\sf\@preprint}
      \end{flushright}
    \end{minipage}
  \end{tabular}
  \vskip 3cm
  \begin{center}
    \LARGE\@title
    \if!\@author!\else \vskip 0.5cm \large\@author\fi
    \if!\@adress!\else \vskip 0.5cm \normalsize\@adress\fi
  \end{center}
  \vskip 2cm
}

\makeatother

\begin{document}

\title{Gerbes and Lie Groups}
\author{Christoph Schweigert and Konrad Waldorf}
\adress{Organisationseinheit Mathematik\\Universit\"at Hamburg\\Bundesstra\ss
e 55\\
D-20146, Germany
}
\preprint{arxiv:0710.5467\\Hamburger Beiträge zur Mathematik Nr. 280\\ZMP-HH/07-10}

\maketitle

\begin{abstract}
Compact Lie groups do not only carry the structure of a
Riemannian manifold, but also canonical families of bundle
gerbes. We discuss the construction of these bundle gerbes and their relation to loop groups.
We present several algebraic structures for bundle gerbes with
connection such as
Jandl structures, gerbe modules and gerbe bimodules, and indicate
their applications to Wess-Zumino terms in two-dimensional
field theories.
\end{abstract}

\section*{Introduction}

Compact Lie groups do not only come with a canonical
metric (the Killing form), but also with a canonical family of bundle gerbes.
These bundle gerbes are geometric objects made of finite dimensional manifolds and maps between
those, and provide a way of 
understanding  structure over the infinite dimensional loop group. 

As a motivation, consider a central extension\index{central extension} of the loop group of a compact connected and simply-connected Lie group
$G$,
\begin{equation*}
0 \longrightarrow  \mathbb{C}^{\times} \longrightarrow  \widehat{LG} \longrightarrow  LG \longrightarrow  0\text{.}
\end{equation*}
Such extensions are classified by  $\mathrm{H}^2(LG,\mathbb{Z})$. 
By transgression, this cohomology group is in turn isomorphic to the cohomology
of $G$,
\begin{equation*}
\mathrm{H}^3(G,\mathbb{Z}) \cong \mathrm{H}^2(LG,\mathbb{Z})\text{.}
\end{equation*}
While the cohomology group $\mathrm{H}^2(LG,\mathbb{Z})$ classifies line
 bundles over $LG$ by their Chern class, $\mathrm{H}^3(G,\mathbb{Z})$ classifies
 bundle gerbes over $G$. In this way, every bundle gerbe over the finite dimensional
 manifold  $G$ gives rise
 to a line bundle over the infinite dimensional manifold $LG$. 
 
\medskip
 
This contribution is  organized as follows. In Section \ref{sec3} we describe
bundle gerbes on general manifolds and their classification. In Section \ref{sec7}
we explain how bundle gerbes can be equipped with a connection
which allows to define surface holonomies: from this point of view,
bundle gerbes generalize the holonomy of principal bundles around curves.  In case that the base
manifold is a compact Lie group $G$, we construct examples of bundle gerbes
over $G$ in Section \ref{sec2}. Then we explain in Section \ref{sec1} how a bundle gerbe gives
rise to a line bundle over the loop space. In Section \ref{sec5} we come
to additional structures for bundle gerbes like bundle gerbe modules,
bimodules and Jandl structures. Finally, in Section \ref{sec6} we outline
the applications of bundle gerbes on Lie groups to two-dimensional conformal field theory and string theory, which are closely related to the surface holonomy
from Section \ref{sec7}. In these theories one also recovers the loop space
as the space of configurations.

\section{Bundle Gerbes}

\label{sec3}

\numberwithin{equation}{section}

Let $M$ be a smooth manifold. 
We shall briefly review  the classification of complex  line bundles over $M$. \index{line bundle}For
this purpose, let us choose a good open cover $\mathfrak{V}=\lbrace V_i \rbrace_{i\in I}$
of $M$, i.e. every finite intersection of open sets $V_i$ is contractible.
In particular, every line bundle $L$ admits  local non-zero sections which
determine   smooth transition functions
\begin{equation}
g_{ij}:V_i \cap V_j \to \mathbb{C}^{\times}\text{.}
\end{equation}
On three-fold intersections $V_i \cap V_j \cap V_k$, these transition functions satisfy the cocycle
condition
\begin{equation}
\label{3}
g_{ik} = g_{ij} \cdot g_{jk}\text{.}
\end{equation}
It is fair to call this equality a cocycle condition, since it means that
$\delta g = 0$ for the element $g:=(g_{ij})\in \check C^1(\mathfrak{V},\mathbb{C}_M^{\times})$
in the \v Cech cohomology of the sheaf of smooth $\mathbb{C}^{\times}$-valued functions on
$M$ with respect to the cover $\mathfrak{V}$. Since we have chosen the cover
$\mathfrak{V}$ to be good, there is a canonical isomorphism $\mathrm{\check H}^1(\mathfrak{V},\mathbb{C}^{\times}_M)\cong
\mathrm{H}^2(M,\mathbb{Z})$ using the exponential sequence. The image of
the class $[g]$ in $\mathrm{H}^2(M,\mathbb{Z})$ is independent of the choice
of the sections, and is called the (first) Chern class $c_1(L)$
of the line bundle $L$. This defines an isomorphism
\begin{equation}
c_1: \mathrm{Pic}(M) \to \mathrm{H}^2(M,\mathbb{Z})
\end{equation}
from the group $\mathrm{Pic}(M)$ of isomorphism classes of line bundles to
the cohomology group $\mathrm{H}^2(M,\mathbb{Z})$, providing a geometric
realization of this group.

\medskip

A bundle gerbe is a geometric object which realizes the cohomology group $\mathrm{H}^3(M,\mathbb{Z})$
in a similar way. To prepare its definition, we fix the
following notation. For a surjective submersion $\pi:Y \to M$ we denote the
$k$-fold fibre product of $Y$ over $M$ by
\begin{equation}
Y^{[k]}= \lbrace (y_1,...,y_k)\in Y^{k}| \pi(y_1)=...=\pi(y_k) \rbrace\text{.}
\end{equation}
This is again a  smooth manifold, having canonical projections
$\pi_{i_1,...,i_\ell}: Y^{[k]} \to Y^{[\ell]}$
on the respective factors. 

\begin{definition}[\normalfont \cite{murray}]
\label{def_gerbe}
A \emph{bundle gerbe}\index{bundle gerbe} $\mathcal{G}$ over
a manifold
$M$ is a triple $(\pi,L,\mu)$ consisting of 
a surjective submersion $\pi: Y
\to M$, a line bundle $L$ over $Y^{[2]}$ and 
an isomorphism
\begin{equation}
\mu : \pi_{12}^{*}L \otimes \pi_{23}^{*}L
\to \pi_{13}^{*}L
\end{equation}
of line bundles over $Y^{[3]}$, such that
$\mu$ is associative
in the sense that the diagram
\begin{equation}
\label{2}
\alxydim{@C=3cm@R=1.5cm}{\pi_{12}^{*}L \otimes \pi_{23}^{*}L \otimes \pi_{34}^{*}L \ar[r]^-{\pi_{123}^{*}\mu\otimes
\id} \ar[d]_{\id \otimes \pi_{234}^{*}\mu} & \pi_{13}^{*}L \otimes \pi_{34}^{*}L
\ar[d]^{\pi_{134}^{*}\mu} \\ \pi_{12}^{*}L \otimes \pi_{24}^{*}L \ar[r]_-{\pi_{124}^{*}
\mu} & \pi_{14}^{*}L}
\end{equation}
of isomorphisms of  line bundles
over $Y^{[4]}$ is commutative. 
\end{definition}

Let us now describe how bundle gerbes realize the cohomology group $\mathrm{H}^3(M,\mathbb{Z})$.
Let us again choose a good open cover $\mathfrak{V}=\lbrace V_i \rbrace_{i\in
I}$ of $M$ which admits sections $s_i: V_i \to Y$ into the manifold $Y$ of
the surjective submersion of a bundle gerbe $\mathcal{G}=(\pi,L,\mu)$. If
we denote by $M_{\mathfrak{V}}$ the disjoint union of all the open sets $V_i$,
the sections $s_i$ patch together to a smooth map $s:M_{\mathfrak{V}} \to
Y$ sending a point $x\in V_i$ to $s_i(x)\in Y$. Note that there are induced
maps $M_{\mathfrak{V}}^{[k]} \to Y^{[k]}$ on fibre products (all denoted
by $s$ in order to simplify the notation) where
$M_{\mathfrak{V}}^{[k]}$ is the disjoint union of all $k$-fold intersections
of open sets $V_i$. Now we pull back the line bundle $L$ along $s$ to a line
bundle over $M_{\mathfrak{V}}^{[2]}$, and the isomorphism $\mu$ to an isomorphism
of line bundles over $M_{\mathfrak{V}}^{[3]}$. For a choice $\sigma_{ij}: V_i \cap V_j \to s^{*}L$ of local sections into the pullback line bundle,
we obtain  smooth functions 
\begin{equation}
g_{ijk}:V_i \cap V_j \cap V_k \to \mathbb{C}^{\times}
\end{equation}
by 
\begin{equation}
s^{*}\mu(\sigma_{ij} \otimes \sigma_{jk}) = g_{ijk} \cdot \sigma_{ik\text{.}}
\end{equation}
The associativity condition (\ref{2}) leads to the equality
\begin{equation}
\label{6}
g_{ijk} \cdot g_{ik\ell} = g_{ij\ell} \cdot g_{jk\ell}
\end{equation}
of functions on four-fold intersections $V_i \cap V_j \cap V_k \cap V_\ell$.
In other words, the element $g=(g_{ijk})\in\check C^2(\mathfrak{V},\mathbb{C}_M^{\times})$
is a cocycle and defines a class in the \v Cech cohomology group $\mathrm{\check H}^2(\mathfrak{V},\mathbb{C}_M^{\times})$.
Its image in the cohomology group $\mathrm{H}^3(M,\mathbb{Z})$ is called
the \emph{Dixmier-Douady class}\index{Dixmier-Douady class} $\mathrm{dd}(\mathcal{G})$ of the bundle gerbe $\mathcal{G}$; it is analogous to the Chern class of a line bundle. 

\medskip

To obtain a classification
result for isomorphism classes of bundle gerbes, we first have to define morphisms between bundle gerbes. To simplify the notation, we work with the
convention that we do not write pullbacks along canonical maps, like in  (\ref{20})
and (\ref{1}) below.        

\begin{definition}
\label{13}
A \emph{morphism}\index{morphism!between bundle gerbes} $\mathcal{A}:\mathcal{G}_1 \to \mathcal{G}_2$ between two bundle gerbes $\mathcal{G}_1=(\pi_1,L_1,\mu_1)$ and $\mathcal{G}_2=(\pi_2,L_2,\mu_2)$ is
a pair $\mathcal{A}=(A,\alpha)$ consisting of  a  vector bundle $A$ over the fibre product $Z:=Y_1 \times_M Y_2$ (whose surjective submersion to $M$
is denoted by $\zeta$) and an isomorphism
\begin{equation}
\label{20}
\alpha: L_1 \otimes \zeta_2^{*}A \to \zeta_1^{*}A \otimes L_2
\end{equation}
of vector bundles over $Z^{[2]}$  such that the  diagram
\begin{equation}
\label{1}
\alxydim{@C=2.5cm}{\zeta_{12}^{*}L_1 \otimes \zeta_{23}^{*}L_1 \otimes \zeta_3^{*}A \ar[r]^-{\mu_1
\otimes \id} \ar[d]_{\id \otimes \zeta_{23}^{*}\alpha} & \zeta_{13}^{*}L_1
\otimes \zeta_3^{*}A \ar[dd]^{\zeta_{13}^{*}\alpha} \\ \zeta_{12}^{*}L_1
\otimes \zeta_2^{*}A \otimes \zeta_{23}^{*}L_2 \ar[d]_{\zeta_{12}^{*}\alpha
\otimes \id} & \\ \zeta_1^{*}A \otimes \zeta_{12}^{*}L_2 \otimes \zeta_{23}^{*}L_2
\ar[r]_-{\id \otimes \mu_2} & \zeta_1^{*}A \otimes \zeta_{13}^{*}L_2}
\end{equation}
of isomorphisms of vector bundles over $Z^{[3]}$ is commutative.
\end{definition}

The definition of the composition of two morphisms $\mathcal{A}:\mathcal{G}_1 \to \mathcal{G}_2$
and $\mathcal{A}':\mathcal{G}_2 \to \mathcal{G}_3$ is quite involved, and
we omit its discussion for the purposes of this article, see \cite{stevenson1,waldorf1}
for more details. Bundle gerbes and
their morphisms fit
into the structure of a 2-category rather than the one of a category. This
becomes obvious when comparing two morphisms $\mathcal{A}$ and $\mathcal{A}'$
between the same bundle gerbes: since $\mathcal{A}$ and $\mathcal{A}'$ consist themselves of vector bundles, the natural way to compare them is a morphism of vector bundles.

\begin{definition}
Let $\mathcal{A}=(A,\alpha)$
and $\mathcal{A}'=(A',\alpha')$
be two morphisms from $\mathcal{G}_1=(\pi_1,L_1,\mu_1)$ to $\mathcal{G}_2=(\pi_2,L_2,\mu_2)$. 
A \index{2-morphism!betw. morph. of bun. ger.}\emph{2-morphism}
\begin{equation}
\beta: \mathcal{A}
\Rightarrow \mathcal{A}'
\end{equation}
is
an
isomorphism $\beta : A \to A'$
of vector bundles over $Z$,
which is compatible with the isomorphisms
$\alpha$ and $\alpha'$ in the sense
that the diagram
\begin{equation}
\alxydim{@C=1.5cm}{
L_1 \otimes \zeta_2^{*}A \ar[r]^-{\alpha}
\ar[d]_{1 \otimes \zeta_2^{*} \beta}
& \zeta_1^{*}A \otimes L_2 \ar[d]^{\zeta_1^{*}\beta\otimes
1}  \\ L_1
\otimes \zeta_2^{*}A' \ar[r]^-{\alpha'}
& \zeta_1^{*}A'  \otimes L_2}
\end{equation}
of isomorphisms of vector bundles
over $Z^{[2]}$ is commutative.
\end{definition}

The 2-categorical setup is also the appropriate context to address the question
which of the morphisms between two bundle gerbes are invertible. 
  
\begin{proposition}[\normalfont \cite{waldorf1}]
\label{14}
A morphism $\mathcal{A}=(A,\alpha)$ is invertible if and only if the vector bundle $A$ is of rank 1. 
\end{proposition}

So, the invertible morphisms in the 2-category are the so-called stable isomorphisms from \cite{murray2}, Section 3.
Let us now return to the relation  between bundle gerbes and the cohomology
group $\mathrm{H}^3(M,\mathbb{Z})$. Let $g$ and $g'$ be the cocycles extracted
from two bundle gerbes $\mathcal{G}$ and $\mathcal{G}'$ over $M$ with respect to the same
cover $\mathfrak{V}$ of $M$ and sections $s_i:V_i \to Y$ and $s_i':V_i
\to Y'$. Now let $\mathcal{A}:\mathcal{G} \to \mathcal{G}'$ be an isomorphism,
$\mathcal{A}=(A,\alpha)$. Let $t:M_{\mathfrak{V}} \to Z:=Y \times_M Y'$
be the map sending a point $x\in V_i$ to the pair $(s_i(x),s_i'(x))$. Using
$t$ we pull back  the line bundle
$A$ and choose non-zero sections $\sigma_i: V_i \to A$. Then we obtain smooth
functions $h_{ij}:V_i \cap V_j \to \mathbb{C}^{\times}$ by
\begin{equation}
\alpha(\sigma_i\otimes \zeta_2^{*}\sigma_{ij}) = h_{ij}\cdot (\zeta_1^{*}\sigma_{ij}
\otimes \sigma_i')\text{.}
\end{equation}
The compatibility condition (\ref{1}) between $\alpha$ and the isomorphisms
$\mu$ and $\mu'$ of the bundle gerbes leads to the equation
\begin{equation}
h_{ik}g_{ijk}=g'_{ijk}h_{ij}h_{jk}\text{,}
\end{equation} 
equivalently, in terms of the \v Cech coboundary operator, $g=g' \cdot \delta
h$. This means that the Dixmier-Douady classes $[g]$ and $[g']$ of the isomorphic bundle
gerbes $\mathcal{G}$ and $\mathcal{G}'$ are equal. This is the main ingredient
for the following classification result.

\begin{theorem}[\normalfont \cite{murray2}]
\label{th1}
The Dixmier-Douady class defines a bijection
between the set of isomorphism
classes of bundle gerbes and the cohomology group $\mathrm{H}^3(M,\mathbb{Z})$.
\end{theorem}

In particular, consider $M=G$ a compact, simple, connected and simply-connected
Lie group. There is a canonical identification $\mathrm{H}^3(G,\mathbb{Z})=\mathbb{Z}$,
so that we have a canonical sequence of bundle gerbes associated to $G$.
In Section \ref{sec2}, we give a geometric construction of these bundle gerbes.

\medskip

The 2-category of bundle gerbes admits several additional structures such
as pullbacks, tensor products
and duals. All these structures are compatible with the Dixmier-Douady class:
\begin{equation}
\mathrm{dd}(\mathcal{G}_1 \otimes \mathcal{G}_2)=\mathrm{dd}(\mathcal{G}_1) + \mathrm{dd}(\mathcal{G}_2)
\quad\text{, }\quad
\mathrm{dd}(\mathcal{G}^{*})=-\mathrm{dd}(\mathcal{G})
\end{equation}
and, for a smooth map $f:X \to M$,
\begin{equation}
\mathrm{dd}(f^{*}\mathcal{G})=f^{*}\mathrm{dd}(\mathcal{G})\text{.}
\end{equation}

To close, let us construct a bundle gerbe whose Dixmier-Douady class vanishes,
representing the neutral element in $\mathrm{H}^3(M,\mathbb{Z})$. For this
purpose, consider the identity $\id_M:M \to M$ as the surjective submersion
and the trivial line bundle $M \times \mathbb{C}$ over $M$. Now, the isomorphism
$\mu$ can be chosen to be the identity $\id_{\mathbb{C}}$,
so that $\mathcal{I}:=(\id_M,M \times \mathbb{C},\id_{\mathbb{C}})$ is a
bundle gerbe. It is easy to verify that its Dixmier-Douady class vanishes,
$\mathrm{dd}(\mathcal{I})=0$. 

\section{Connections on Bundle Gerbes and Holonomy}

\label{sec7}

Before we discuss more examples of bundle gerbes in Section \ref{sec2} we introduce
several additional structures on bundle gerbes and the appropriate cohomology
theory for their classification.

Again, we first review similar additional structures for complex line bundles. One can
equip every line bundle with a hermitian metric and a (unitary)
connection $\nabla$.\index{connection!on a line bundle} Additionally to the transition function $g_{ij}:V_i \cap V_j
\to \mathbb{C}^{\times}$, which now can be determined such that it takes values
in $U(1)$, we obtain local connection 1-forms $A_i\in\Omega(V_i)$ by writing
the covariant derivatives of the sections $s_i:V_i \to L$ as
$\nabla s_i = A_{i} \otimes s_i$. 
The Leibniz rule implies the equality
\begin{equation}
\label{4}
A_j - A_i - \mathrm{dlog}(g_{ij})=0\text{.}
\end{equation}
As we shall see next, the local data $(g,A)$ define a cocycle in the Deligne complex $\mathcal{D}_{\mathfrak{V}}^k(n)$ for $n=1$.
\index{Deligne cohomology}The first  cochain groups of this complex are 
\begin{eqnarray}
\label{9}
\mathcal{D}_{\mathfrak{V}}^0(1)&=&\check C^0(\mathfrak{V},U(1)_M) 
\\
\mathcal{D}_{\mathfrak{V}}^1(1)&=&\check C^1(\mathfrak{V},U(1)_M) \oplus \check C^0(\mathfrak{V},\Omega^1_M)
\\
\mathcal{D}_{\mathfrak{V}}^2(1)&=&\check C^2(\mathfrak{V},U(1)_M) \oplus \check C^1(\mathfrak{V},\Omega^1_M)\text{,}
\end{eqnarray}
and its differential is given
by
\begin{equation}
\mathrm{D}: \mathcal{D}^1_{\mathfrak{V}}(1) \to \mathcal{D}^2_{\mathfrak{V}}(1):(g,A)
\mapsto (\delta g, \delta A - \mathrm{dlog}(g))\text{.}
\end{equation}
In general, the Deligne complex $\mathcal{D}_{\mathfrak{V}}^k(n)$ is the total complex of the \v Cech-Deligne double complex, whose rows are \v Cech cochain
groups of the sheaf complex
\begin{equation}
\alxy{0 \ar[r] & U(1)_M \ar[r]^-{\mathrm{dlog}} & \Omega^1_M \ar[r]^-{\mathrm{d}} &... \ar[r]^-{\mathrm{d}} & \Omega^n_M\text{,}}
\end{equation}
truncated in degree $n$, and the columns are the usual \v Cech complexes
associated to the sheaves $U(1)_M$ and $\Omega_M^k$. The truncation is necessary to
describe not just flat objects.

So we can regard the local data $(g,A)$ of a hermitian line bundle\index{line
bundle} with  connection as an element of $\mathcal{D}_{\mathfrak{V}}^1(1)$. By equations (\ref{3}) and (\ref{4}) it satisfies $\mathrm{D}(g,A)=0$, and thus defines a class in
the first cohomology group of the Deligne complex, denoted by $\mathrm{H}^1(M,\mathcal{D}(1))$.
One can show that this  defines a bijection
\begin{equation}
\label{12}
\mathrm{Pic}^{\nabla}(M) \cong \mathrm{H}^1(M,\mathcal{D}(1))
\end{equation}
between the group $\mathrm{Pic}^{\nabla}(M)$ of isomorphism classes of hermitian
line bundles with connection  and Deligne cohomology \cite{brylinski1}. 

\medskip
 
Now we discuss bundle gerbes with similar additional structures. 
\begin{definition}
Let $\mathcal{G}=(\pi,L,\mu)$ be a bundle gerbe over $M$.
It can be equipped successively with the following additional structures.\begin{enumerate}
\item[(a)]
A \emph{hermitian structure} on $\mathcal{G}$ is a hermitian metric on the line
bundle $L$, such that the isomorphism $\mu$ is an isometry.
\item[(b)]
A \emph{connection}\index{connection!on a bundle gerbe} on the hermitian bundle gerbe $\mathcal{G}$ is a connection
$\nabla$
on the hermitian line bundle $L$, such that the isomorphism $\mu$ respects the induced
connections. 
\item[(c)]
A \emph{curving}\index{curving} of a connection $\nabla$ on the hermitian bundle gerbe $\mathcal{G}$ is a 2-form $C\in\Omega^2(Y)$, such that
\begin{equation}
\label{11}
\pi_2^{*}C-\pi_1^{*}C=\mathrm{curv}(\nabla)\text{,}
\end{equation}
where $\mathrm{curv}(\nabla) \in \Omega^2(Y^{[2]})$ is the curvature of the
connection $\nabla$ on $L$.
\end{enumerate}
\end{definition}

Every bundle gerbe admits all of these additional structures \cite{murray}.
Because the applications of bundle gerbes we discuss later require all this
additional structure, we work from now only with hermitian bundle
gerbes with connection and curving. 

An important feature of those gerbes
is that they provide a notion of curvature. To see this, consider the derivative
of equation (\ref{11}): since the curvature of the connection $\nabla$ is a closed
form, we obtain $\pi_1^{*}\mathrm{d}C=\pi_2^{*}\mathrm{d}C$, which means
that $\mathrm{d}C$ is the pullback of a 3-form on $M$, 
\begin{equation}
\mathrm{d}C=\pi^{*}H\text{.}
\end{equation}
This 3-form $H$ is uniquely determined and closed; it is called the curvature
of the curving of the connection $\nabla$ on the hermitian  bundle gerbe $\mathcal{G}$,
and denoted by $\mathrm{curv}(C):=H$.

To have a simple example of such additional structures, the bundle gerbe $\mathcal{I}=(\id_M,M \times \mathbb{C},\id_{\mathbb{C}})$ with vanishing Dixmier-Douady
class becomes a hermitian bundle gerbe with connection by taking the canonical
hermitian metric and the trivial flat connection $\nabla:=\mathrm{d}$ on the trivial line bundle $M\times
\mathbb{C}$. Note that now any 2-form $C\in\Omega^2(M)$ satisfies the condition
(\ref{11}) for a curving on $I$, because $\nabla$ is flat and $\pi_1=\pi_2=\id_{M}$.
So we have a canonical hermitian bundle gerbe $\mathcal{I}_C:=(\id_M,M \times
\mathbb{C},\id_{\mathbb{C}},\mathrm{d},C)$ with connection 
and curving for every 2-form $C \in \Omega^2(M)$. The curvature of its curving
$C$ is $\mathrm{curv}(C)=\mathrm{d}C$.

\medskip

Now we extend the cohomological classification from bundle gerbes to hermitian
bundle gerbes  $\mathcal{G}=(\pi,L,\mu,\nabla,C)$  with connection and curving, using a  good open cover $\mathfrak{V}$ of $M$. Recall that we have extracted
a transition function $g_{ijk}:V_{i}\cap V_j \cap V_k \to \mathbb{C}^{\times}$
using a choice of sections $s_i:V_i \to Y$ and $\sigma_i: V_i \to s^{*}L$,
defining a representative $g\in \check C^2(\mathfrak{V},\mathbb{C}^{\times})$
of the Dixmier-Douady class of the bundle gerbe. Now that $L$ is a hermitian
line bundle, we choose the sections $\sigma_i$ such that $g_{ijk}$ is $U(1)$-valued.
Furthermore, by using the connection $s^{*}\nabla$ on $s^{*}L$, we obtain
local connection 1-forms $A_{ij}\in \Omega^1(V_i \cap V_j)$. The condition
that $\mu$ preserves connections implies 
\begin{equation}
\label{7}
A_{jk}-A_{ik}+A_{ij}+\mathrm{dlog}(g_{ijk})=0\text{.}
\end{equation}
Finally, the
curving $C$ gives rise to local 2-forms $B_i:= s_i^{*}C\in\Omega^2(V_i)$,
and
the compatibility (\ref{11}) with the curvature of $\nabla$ implies
\begin{equation}
\label{8}
B_j - B_i - \mathrm{d}A_{ij}=0\text{.}
\end{equation}
Note that the curvature $\mathrm{curv}(C)$ of the curving can be computed from the local data by $H|_{V_i}=\mathrm{d}B_i$.
By (\ref{8}), this gives indeed a globally defined 3-form.

In terms of \v Cech cohomology, we have extracted an element $(g,A,B)$ in
the second  cochain group of the Deligne complex in degree 2,
\index{Deligne cohomology}
\begin{equation}
\label{10}
\mathcal{D}^2_{\mathfrak{V}}(2) = \check C^2(\mathfrak{V},U(1)) \oplus \check C^1(\mathfrak{V},\Omega^1)
\oplus \check C^0(\mathfrak{V},\Omega^2)\text{.}
\end{equation}
The differential is here
\begin{equation}
\mathrm{D}:\mathcal{D}_{\mathfrak{V}}^2(2) \to \mathcal{D}^3_{\mathfrak{V}}(2):(g,A,B) \mapsto (\delta
g, \delta A + \mathrm{dlog}(g),\delta B-\mathrm{d}A)\text{,}
\end{equation}
so that the cocycle condition (\ref{11}) and equations (\ref{7}) and (\ref{8})
imply the cocycle condition $\mathrm{D}(g,A,B)=0$. This way, a hermitian bundle
gerbe with connection and curving defines a class in the cohomology of the
Deligne complex in degree 2, $\mathrm{H}^2(M,\mathcal{D}(2))$. As an exercise,
the reader may compute the Deligne class of the canonical hermitian bundle
gerbe $\mathcal{I}_C$ with connection and curving associated to any 2-form
$C\in\Omega^2(M)$.

Note that both the Deligne cochain groups $\mathcal{D}_{\mathfrak{V}}^1(1)$ from (\ref{9}) and
$\mathcal{D}_{\mathfrak{V}}^2(2)$ from (\ref{10}) have projections on the first summand,
which commute with the Deligne differential and the \v Cech coboundary operator,
so
that we get induced (surjective) group homomorphisms \cite{brylinski1}
\begin{equation}
\mathrm{H}^k(M,\mathcal{D}(k)) \to \mathrm{H}^{k+1}(M,\mathbb{Z})\text{.}
\end{equation}
This way we obtain the Chern class and the Dixmier-Douady class of a hermitian
line bundle with connection and of a hermitian bundle gerbe with connection
and curving,  respectively, from their Deligne classes. Its surjectivity means that Deligne
cohomology refines the ordinary singular cohomology with $\mathbb{Z}$ coefficients.

To achieve a classification result for hermitian bundle gerbes with connection
and curving similar to the result (\ref{12}) for hermitian line bundles with
connection, we have to adapt the definition of a morphism between bundle
gerbes to morphisms between hermitian bundle gerbes with connection and curving.

\begin{definition}
A \emph{morphism}\index{morphism!between bundle gerbes} $\mathcal{A}:\mathcal{G} \to \mathcal{G}'$ between two hermitian
bundle gerbes $\mathcal{G}=(\pi,L,\mu,\nabla,C)$ and $\mathcal{G}'=(\pi',L',\mu',\nabla',C')$
with connection and curving is a morphism $\mathcal{A}=(A,\alpha)$
as in Definition \ref{13}, together with a connection $\blacktriangledown$
on the vector bundle $A$, such that
\begin{enumerate}
\item 
the isomorphism $\alpha$ respects connections.
\item
 the curvature of $\blacktriangledown$ is related to the curvings
by
\begin{equation}
\label{23}
\frac{1}{n}\mathrm{tr}(\mathrm{curv}(\blacktriangledown))= C' - C\text{.}
\end{equation} 
\end{enumerate}
\end{definition}

As in Proposition \ref{14}, a morphism is invertible precisely if the vector bundle
is of rank 1. It is again an easy exercise to check, that an isomorphism $\mathcal{A}:\mathcal{G}
\to \mathcal{G}'$ of hermitian bundle gerbes with connection and curving
with Deligne cocycles $(g,A,B)$ and $(g',A',B')$ respectively, gives rise to
a Deligne cochain $(h,W)\in \mathcal{D}_{\mathfrak{V}}^1(2)$ which satisfies
\begin{equation}
(g',A',B') = (g,A,B) + \mathrm{D}(h,W)\text{.}
\end{equation}
This shows that isomorphism classes of  hermitian bundle gerbes with connection
and curving  have well-defined Deligne classes.

\begin{theorem}[\normalfont \cite{murray2}]
Isomorphism classes of hermitian bundle gerbes with connection and curving
correspond bijectively to the Deligne cohomology group $\mathrm{H}^2(M,\mathcal{D}(2))$.
\index{Deligne cohomology}
\end{theorem}

Particular examples of morphisms are trivializations: a \emph{trivialization}\index{trivialization!of
a bundle gerbe} of
a hermitian bundle gerbe $\mathcal{G}$ with connection and curving is an isomorphism
\begin{equation}
\mathcal{T}:\mathcal{G} \to \mathcal{I}_{\rho}
\end{equation}
for some 2-form $\rho\in\Omega^2(M)$. In terms of local data, this
isomorphism corresponds to a Deligne cochain $(h,W)\in\mathcal{D}_{\mathfrak{V}}^1(2)$ with
\begin{equation}
\label{16}
(1,0,\rho) = (g,A,B) + \mathrm{D}(h,W)\text{,}
\end{equation}
if $(g,A,B)$ is local data of the bundle gerbe $\mathcal{G}$. In particular,
the existence of a trivialization implies that the Dixmier-Douady class of
$\mathcal{G}$ vanishes. Many assertions
about bundle gerbes and their isomorphisms
can be proven either in a geometrical way or by computations in Deligne cohomology.
As an illustration, we shall prove the following important 
\begin{lemma}
\label{lem1}
Two trivializations 
\begin{equation}
\mathcal{T}_1:
\mathcal{G} \to \mathcal{I}_{\rho_1}
\quad\text{and}\quad
\mathcal{T}_2: \mathcal{G}
\to \mathcal{I}_{\rho_2}
\end{equation}
 of the
same hermitian bundle gerbe $\mathcal{G}$
over $M$ with connection and curving determine
a hermitian line bundle over $M$ with connection of curvature
$\rho_2-\rho_1$.
\end{lemma}

\textit{Proof 1 (2-categorical).}
Using the features of the 2-category of bundle
gerbes, we can give a  quite conceptual proof:
by taking the inverse and composition (which we have not explained in this
article, but can be found in  \cite{waldorf1}), we
obtain an isomorphism
\begin{equation}
\mathcal{T}_2 \circ \mathcal{T}_1^{-1}: \mathcal{I}_{\rho_1}
\to \mathcal{I}_{\rho_2}
\end{equation}
of trivial bundle gerbes. From the definitions
of isomorphisms and trivial bundle gerbes
it follows immediately, that $\mathcal{T}_2 \circ \mathcal{T}_1^{-1}$ is a line bundle
with curvature $\rho_2-\rho_1$.
\endofproof

\textit{Proof 2 (geometrical).}
The two isomorphisms  $\mathcal{T}_i=(T_i,\tau_i,\blacktriangledown_i)$ consist
of hermitian line bundles
$T_i$ over $Z:=Y \times_M M \cong Y$, connections $\blacktriangledown_i$
of curvature $\mathrm{curv}(\blacktriangledown_i)=\pi^{*}\rho_i - C$, and isomorphisms $\tau_i: L \otimes
\pi_2^{*}T_i \to \pi_1^{*}T_i$ of hermitian line bundles
respecting the connections. They can be
composed to an isomorphism
\begin{equation}
\tau_2^{-1} \otimes \tau_1^{*}: \pi_1^{*}(T_2 \otimes T_1^{*}) \to \pi_2^{*}(T_2
\otimes T_1^{*})
\end{equation}
of hermitian line bundles with connection over $Y^{[2]}$, which satisfies
the obvious cocycle condition over $Y^{[3]}$, due to the commutative diagram
(\ref{1}) for the $\tau_i$. This is the condition for the hermitian line bundle $T_2 \otimes T_1^{*}$ with connection $\blacktriangledown_2 - \blacktriangledown_1$
to descend from $Y$ to $M$. The descent line bundle
has the claimed curvature. 
\endofproof

\textit{Proof 3 (cohomological).} If the isomorphisms $\mathcal{T}_1$ and $\mathcal{T}_2$ have local data $(h_1,W_1)$ and
$(h_2,W_2)$ respectively, both satisfying equation (\ref{16}), their difference
satisfies
\begin{equation}
\mathrm{D}(h_2 \cdot h_1^{-1},W_2 - W_1) = 0\text{,}
\end{equation}
which is the cocycle condition for a hermitian line bundle with connection
of curvature $\mathrm{d}(W_2-W_1)=\mathrm{curv}(\blacktriangledown_2)-\mathrm{curv}(\blacktriangledown_1)=\rho_2
- \rho_1$. 
\endofproof

One of the most important aspects of the theory of hermitian bundle gerbes
with connection and curving is that
they provide a notion of holonomy around surfaces. 

\begin{definition}[\normalfont \cite{carey2}]
\label{def2}
Let $\mathcal{G}$ be a hermitian bundle gerbe with connection and curving over
$M$. For a closed oriented surface $\Sigma$
and a smooth map $\phi: \Sigma
\to M$, let  
\begin{equation}
\mathcal{T}:\phi^{*}\mathcal{G}
\to \mathcal{I}_{\rho}
\end{equation}
be a trivialization of the pullback of 
$\mathcal{G}$ along $\phi$.
Then we define
\begin{equation}
\mathrm{hol}_{\mathcal{G}}(\phi)
:= \exp \left (\int_{\Sigma} \rho
\right
)
\end{equation}
to be the \emph{holonomy}\index{holonomy!of a bundle gerbe} of  $\mathcal{G}$ around $\phi:\Sigma
\to M$. 
\end{definition}

Note that the Dixmier-Douady class of $\phi^{*}\mathcal{G}$ vanishes by dimensional
reasons, so that the existence of the trivialization $\mathcal{T}$ is guaranteed.
However, it is not unique, and different trivializations
may have different 2-forms $\rho$. Now, by Lemma \ref{lem1} we know that
the difference $\rho_2-\rho_1$
between two such 2-forms is the
curvature of some line bundle over
$M$,
in particular: it is a closed form
with  integral class. Then, the calculation
\begin{equation}
\exp \left ( \int_\Sigma \rho_2 \right
) = \exp \left ( \int_\Sigma \rho_2-\rho_1 \right
) \cdot \exp \left ( \int_\Sigma \rho_1 \right
)= \exp \left ( \int_\Sigma \rho_1 \right
)
\end{equation}  
shows that the definition $\mathrm{hol}_{\mathcal{G}}(\phi)$
is independent of the choice of
the trivialization. 

There also exist expressions for the holonomy $\mathrm{hol}_{\mathcal{G}}(\phi)$
in terms of local data of $\mathcal{G}$. They generalize the local formulae
for the holonomy of hermitian line bundles with connection.   In Section \ref{sec6} we describe
the applications of holonomy of bundle gerbes in conformal field theory.

\section{Bundle Gerbes over compact Lie Groups}

\label{sec2}

Now that we have introduced bundle gerbes as geometric objects over arbitrary
manifolds,  we specialize to manifolds which are Lie groups. We describe
in this section, how the Lie group structure allows constructions of examples of bundle gerbes. First constructions of  gerbes over different types of compact Lie groups
(in  realizations different from bundle gerbes) can be found in \cite{chatterjee,brylinski2}.
Bundle gerbes (with connection and curving)\index{bundle gerbe} have been constructed in \cite{gawedzki1,meinrenken1,gawedzki2}.

\medskip

As already mentioned before, for a compact, simple, connected and simply-connected Lie group $G$ \index{Lie group!simply-connected} have $\mathrm{H}^3(G,\mathbb{Z})=\mathbb{Z}$.
The (up to isomorphism unique) bundle gerbe over $G$ whose Dixmier-Douady class corresponds
to $1\in \mathbb{Z}$ is called the basic bundle gerbe, \index{bundle gerbe!basic} and denoted by $\mathcal{G}_0$.
The bundle gerbe with Dixmier-Douady class $k\in \mathbb{Z}$ can then be
obtained from $\mathcal{G}_0$ or $\mathcal{G}_0^{*}$ by a $k$-fold tensor
product.

For the purposes of this article, we will restrict ourselves to the  construction
given in \cite{gawedzki1}, by which we obtain the basic bundle gerbe over
the special unitary groups
$\mathrm{SU}(n)$ and the symplectic groups $\mathrm{Sp}(n)$. 
First we consider a general compact, simple and simply-connected Lie group $G$ with
Lie algebra $\mathfrak{g}$. We choose a maximal torus $T$ with Lie algebra
$\mathfrak{t}$ of rank $r$. 
 We further choose a set of simple roots $\alpha_1,...,\alpha_r$ and denote
 the associated positive Weyl chamber by $C$. Let $\alpha_0$ the highest
 root and let 
\begin{equation}
\mathfrak{A}:= \lbrace \xi\in C\;|\;\alpha_0(\xi)\leq1 \rbrace\text{.}
\end{equation}
be the fundamental alcove. It is  bounded by the hyperplanes
$H_i$ perpendicular to the roots $\alpha_i$, and the additional hyperplane $H_0$
consisting of elements $\xi\in\mathfrak{t}$ with $\alpha_0(\xi)=1$. So it is
an
$r$-dimensional simplex
in $\mathfrak{t}$ with vertices
$\mu_0...,\mu_{r}$ determined
by the condition that $\mu_i\in H_{j}$ for
all $j \neq i$.  

For simple and
simply connected groups, the fundamental alcove parameterizes conjugacy
classes of $G$ in the sense that
each conjugacy class contains a
unique point $\exp\xi$
with $\xi\in\mathfrak{A}$. This
defines a continuous map 
\begin{equation}
q:G \to \mathfrak{A}\text{.}
\end{equation}
Let
$\mathfrak{A}_i$ be the open complement
of the face opposite to the vertex
$\mu_i$ in $\mathfrak{A}$, and
consider the open sets $U_i :=
q^{-1}(\mathfrak{A}_i)$. More generally, for any subset $I\subset \underline{r}=\lbrace
0,...,r\rbrace$ denote by $U_{I}$
the intersection of all $U_i$ with
$i\in I$, and similarly by $\mathfrak{A}_{I}$
the intersection of all $\mathfrak{A}_{i}$
with $i\in I$. Of course $U_{I}=q^{-1}(\mathfrak{A}_I)$.
We use the open sets $U_i$ to construct
the surjective submersion
\begin{equation}
Y:= \bigsqcup_{i\in \underline{r}}
U_i
\quad\text{ and }
\quad
\pi:Y \to G: (i,x) \mapsto x\text{.}
\end{equation}
Note that the $k$-fold fibre products
are disjoint unions of intersections
\begin{equation}
Y^{[k]} = \bigsqcup_{|I|=k} U_I\text{.}
\end{equation}
The surjective submersion $\pi:Y \to M$ will serve as the first ingredient
of the bundle gerbe we want to construct. 
To construct the line bundle $L$ over $Y^{[2]}$ we  show  next that  $Y^{[2]}$
projects onto a  union of coadjoint orbits. 

For any $I \subset \underline{r}$, all group elements $\exp\xi$ with $\xi$
in the open face spanned by the vertices $\mu_i$ with $i\in I$ have the same
centralizer $G_I$. 
For any inclusion $I\subset
J$ it follows  that $G_J \subset G_I$;
for $I=\underline{r}$ we obtain
$G_{\underline{r}}=T$. Let $S_I$ be the orbit of $\exp\mathfrak{A}_{I} \subset T$ under the conjugation with
$G_{I}$. Consider the set $G \times_{G_I} S_I$ consisting of equivalence
classes of pairs $(g,s)\in G \times S_I$ under the equivalence relation $(g,s)\sim
(gh,h^{-1}sh)$ for $h\in G_I$. We have the canonical  projection $\rho_{I}:
G \times_{G_I} S_I \to G/G_{I}$ and a smooth map
\begin{equation}
u_I: G \times_{G_I} S_I \to U_I
\end{equation}
which sends a representative $(g,s)$ to $gsg^{-1}\in G$. This is well-defined
on equivalence classes, and for $h\in G_I$ and $\xi\in \mathfrak{A}_I$ with $s:=h\exp\xi \,h^{-1}$ we find $q(gsg^{-1})=\xi$ and hence $gsg^{-1}\in
U_I$. The map $u_I$ is even a diffeomorphism: for $g\in U_I$ let $\xi:=q(g)\in\mathfrak{A}_I$
and $h\in G$ such that $g=h\exp\xi \,h^{-1}$. Then, the inverse sends $g$ to
the equivalence class of $(h,\exp \xi)$.

Since $G_{ij}$ fixes the difference $\mu_{ij}:=\mu_j-\mu_i$, the quotient $G/G_{ij}$
projects on the coadjoint orbit $\mathcal{O}_{ij}$
of $\mu_{ij}$ in $\mathfrak{g}^{*}$. \index{coadjoint orbit}
Now we specialize our construction to the case that $\mu_{ij}$ is a weight.
This is the case for  $\mathrm{SU}(n)$ and  $\mathrm{Sp}(n)$, where the vertices of the fundamental alcove
are contained in the weight lattice. For $\mu_{ij}$ being a weight there is a canonical line bundle
$L_{ij}$ over the coadjoint orbit $\mathcal{O}_{{ij}}$. \index{line bundle!over
a coadjoint orbit}Let us recall
the construction of the pullback of this line bundle to $G/G_{ij}$.
Let $\chi_{ij}: G_{ij} \to U(1)$
be the character associated to
the weight $\mu_{ij}$. Then the
line bundle is the bundle associated
to the principal $G_{ij}$-bundle
$G$ over $G/G_{ij}$, namely
\begin{equation}
\label{b1}
L_{ij} := G \times_{G_{ij}}
\mathbb{C}\text{.}
\end{equation}
The line
bundles $L_{ij}$ over
$G/G_{ij}$ can be pulled
back along
\begin{equation}
\alxy{
U_i \cap U_j \ar[r]^-{\sim} &  G \times_{G_{ij}}
S_{ij} \ar[r] & G/G_{ij}}
\end{equation}
to line bundles over $U_i \cap
U_j$, and their disjoint union
gives a line bundle $L$ over $Y^{[2]}$.

To close the definition of a bundle gerbe over $G$, it remains to  construct
 the isomorphism $\mu$ of line
bundles over $Y^{[3]}$, i.e. we need an isomorphism
\begin{equation}
\mu: \pi_{12}^{*}L \otimes \pi_{23}^{*}L \to \pi_{13}^{*}L
\end{equation}
of line bundles over $Y^{[3]}$. Over each connected component $U_i \cap U_j
\cap U_k$, this can be chosen as the pullback of the canonical identification
\begin{equation}
L_{ij} \otimes L_{jk} \cong L_{ik}
\end{equation}
of line bundles over $G/G_{ijk}$, which comes from the coincidence $\mu_{ik}=\mu_{ij}+\mu_{jk}$ of the weights from which the line bundles are constructed. This identification  is obviously
associative.

\medskip

In order to calculate the  Dixmier-Douady class of the bundle gerbe we just
have constructed,
we choose a connection and a curving. This procedure is analogous to
the calculation of Chern classes of complex vector bundles using a choice
of a hermitian
metric and a connection, see \cite{milnor1}. Recall that this is based on
the fact that the de Rham cohomology class of the curvature of a hermitian
line bundle with connection equals the image of its Chern class under the induced
map $\iota^{*}:\mathrm{H}^2(M,\mathbb{Z}) \to \mathrm{H}^2(M,\mathbb{R})$. The same
holds for a bundle gerbe $\mathcal{G}=(\pi,L,\mu,\nabla,C)$:
\begin{equation}
[\mathrm{curv}(C)] = \iota^{*}\mathrm{dd}(\mathcal{G})\text{,}
\end{equation}
which can be proven by a zigzag-argument in the \v Cech-Deligne double complex
\cite{brylinski1}. 

Note that -- in a certain normalization of the Killing form \cite{pressley1} -- the  bi-invariant closed 3-form
 \begin{equation}
H:=\frac{1}{6} \left \langle \theta \wedge [\theta \wedge \theta]  \right \rangle \in
\Omega^3(G)
\end{equation} 
is a generator of the cohomology group $\mathrm{H}^3(G,\mathbb{Z})$. Here $\theta$
is the left invariant Maurer-Cartan form, a $\mathfrak{g}$-valued 1-form  on $G$; the 3-form $H$ takes at the unit element the value $H_1(X,Y,Z) = \left \langle X,[Y,Z]  \right \rangle$ on elements $X,Y,Z\in \mathfrak{g}$.
Our goal
is now to define a connection and a curving with curvature
$H$.

First note that the line bundle
$L_{ij}$ from
(\ref{b1}) inherits a hermitian metric from the standard metric on $\mathbb{C}$,
and that the isomorphism $\mu$ is an isometry. The line bundle $L_{ij}$ can
also be equipped with a connection: we consider
 the 1-form $A_{ij}:=\left \langle \mu_{ij},\theta  \right \rangle$ on the
 total space of the  principal $G_{ij}$-bundle $G$. It induces a connection on
 the associated line bundle
 $L_{ij}$ because   $\mu_{ij}$
is preserved under the action of $G_{ij}$. This way the  bundle gerbe is a hermitian bundle gerbe with connection. 

To define a curving for this bundle gerbe, i.e. a 2-form $C \in \Omega^2(Y)$,
we use the fact that the linear retraction of $\mathfrak{A}_i$ to the vertex $\mu_i$
lifts to a smooth retraction of $U_i$ to the conjugacy class $\mathcal{C}_{\mu_i}$
\cite{meinrenken1}, \index{conjugacy class}
\begin{equation}
r_i: U_i \times [0,1] \to U_i\text{.}
\end{equation}
On the conjugacy class $\mathcal{C}_{\mu_i}$, the 3-form $H$ becomes exact,
$\iota_i^{*}H = \mathrm{d}\omega_{\mu_i}$ where $\iota_i: \mathcal{C}_{\mu_i}
\to G$ is the inclusion, and 
\begin{equation}
  \omega_{\mu_i} := \Big \langle \iota_i^{*}\theta \wedge \frac{\mathrm{Ad}^{-1}
  + \id_{\mathfrak{g}}}{\mathrm{Ad}^{-1}-\id_{\mathfrak{g}}}\iota_i^{*}\theta\Big \rangle
  \in \Omega^2(\mathcal{C}_{\mu_i})
\end{equation}
is an invariant 2-form on the conjugacy class \cite{bordalo1}. Here the notation is to be understood as follows: at some element $h\in \mathcal{C}_{\mu_i}$, the 2-form $\omega_{\mu_i}$ is obtained by
considering the Maurer-Cartan forms in $h$ and taking the inverse of the adjoint
action $(\mathrm{Ad}_h)^{-1}: \mathfrak{g} \to \mathfrak{g}$. The endomorphism $(\mathrm{Ad}_h)^{-1}-\id_{\mathfrak{g}}$ becomes invertible when restricted to the image of $\iota_i^{*}\theta$ so that the fraction makes sense.

By pullback along $r_i$ and fibre integration, one obtains a 2-form $C_i
\in \Omega^2(U_i)$ with \begin{equation}
H|_{U_i}=\mathrm{d}C_i \text{.}
\end{equation}
One can now show that $C_j - C_i = \left \langle \mu_{ij},\mathrm{d}\theta  \right \rangle$ on $U_i \cap U_j$, which is the condition on the curving.
 So, our construction realizes, for $G=\mathrm{SU}(n)$
and $G=\mathrm{Sp}(n)$ a hermitian bundle gerbe with connection and curving
with Dixmier-Douady class $1 \in \mathbb{Z}$, the basic bundle gerbe.

For the other compact, simple, connected and
simply-connected Lie groups, one can find an integer $k_0$ for which the
vertices of $k_0\mathfrak{A}$ are weights. Using the weights $k_0\mu_{ij}$
in the construction of the line bundles $L_{ij}$, and the 2-forms $k_0C_i$
in the definition of the curving, we obtain bundle gerbes
with Dixmier-Douady class  $k_0\in \mathbb{Z}$ \cite{meinrenken1}. The smallest such integer $k_0$ is tabulated in \cite{bourbaki1}:
\begin{center}\begin{tabular}{|c|c|c|c|c|c|c|c|c|c|}\hline
$G$ & $\mathrm{SU(n)}$ & $\mathrm{Spin}(n)$ & $\mathrm{Sp}(2n)$ & $E_6$ & $E_7$ & $E_8$ & $F_4$ & $G_2$  \\ \hline
$k_0$ & 1 & 2 & 1 & 3 & 12 & 60 & 6 & 2 \\ \hline
\end{tabular}\end{center}
The construction of  the basic bundle gerbes on the groups with $k_0>1$ requires more advanced techniques \cite{meinrenken1,gawedzki2}. Here it becomes
in particular important that the definition of a bundle gerbes admits $\pi:Y
\to G$ to be a surjective submersion, which is more general than just an
open cover of $G$.

\medskip

Starting from a bundle gerbe $\mathcal{G}$ on a simply-connected Lie group $G$, we now describe a method to obtain bundle gerbes on the non simply-connected
Lie groups $G/Z$ which are quotients of $G$ by a subgroup $Z$ of the center
of $G$. More generally, let $\Gamma$ be a finite group acting on a manifold $M$ by
diffeomorphisms.

\begin{definition}
\label{def1}
A \emph{$\Gamma$-equivariant structure}\index{equivariant structure!on a
bundle gerbe} $(\mathcal{A}_{\gamma},\varphi_{\gamma_1,\gamma_2})$ on a bundle gerbe   $\mathcal{G}$ over $M$ consists of
isomorphisms
\begin{equation}
\mathcal{A}_{\gamma}: \mathcal{G}
\to \gamma^{*}\mathcal{G}
\end{equation}
for each $\gamma\in\Gamma$ and of 2-isomorphisms
\begin{equation}
\varphi_{\gamma_1,\gamma_2}: \gamma_1^{*}\mathcal{A}_{\gamma_2}
\circ \mathcal{A}_{\gamma_1}
\Longrightarrow \mathcal{A}_{\gamma_2\gamma_1}
\end{equation}
for each pair $\gamma_1,\gamma_2\in
Z$, such that the diagram
\begin{equation}
\label{5}
\alxydim{@C=2.4cm@R=1.5cm}{\gamma_1^{*}\gamma_2^{*}\mathcal{A}_{\gamma_3} \circ
\gamma_1^{*}\mathcal{A}_{\gamma_2}
\circ \mathcal{A}_{\gamma_1} \ar@{=>}[d]_{\gamma_1^{*}\varphi_{\gamma_2,\gamma_3}
 \circ \id} \ar@{=>}[r]^-{\id \circ \varphi_{\gamma_1,\gamma_2}}
& \gamma_1^{*}\gamma_2^{*}\mathcal{A}_{\gamma_3}
\circ \mathcal{A}_{\gamma_2\gamma_1} \ar@{=>}[d]^{\varphi_{\gamma_2\gamma_1,\gamma_3}}
\\ \gamma_1^{*}\mathcal{A}_{\gamma_3\gamma_2} \circ
\mathcal{A}_{\gamma_1}
\ar@{=>}[r]_-{\varphi_{\gamma_1,\gamma_3\gamma_2}}
& \mathcal{A}_{\gamma_3\gamma_2\gamma_1}}
\end{equation}
of 2-isomorphisms is commutative. 
\end{definition}

Not every bundle gerbe $\mathcal{G}$ over $M$ admits a $\Gamma$-equivariant structure,
and if it does, it may not be unique. To obtain obstructions and classifications
we use the cohomological language. For this purpose, we impose the structure of a $\Gamma$-module on the Deligne cochain groups $\mathcal{D}_{\mathfrak{V}}^n(k)$ \cite{gawedzki6}. We assume the existence of a
good open cover $\mathfrak{V}=\lbrace V_i \rbrace_{i\in I}$ of $M$ which is compatible with the group
action in the sense that there is an induced action of $\Gamma$ on the index set
$I$ such that $\gamma(V_i)=V_{\gamma i}$. For example, the open sets $U_i$ we have used
in the construction of the basic bundle gerbe satisfy this condition. Then $\Gamma$ acts by pullback on the cochain
groups $\mathcal{D}_{\mathfrak{V}}^n(k)$. 

For each $\Gamma$-module $W$, one can build the usual group cohomology
 complex, consisting of
cochain groups $C^p_{\Gamma}(W) := \mathrm{Map}(\Gamma^{p+1},W)$ and the usual coboundary operator \index{group cohomology}
\begin{equation}
\mathrm{d}: C_{\Gamma}^{k-1}(W) \to C_{\Gamma}^{k}(W)\text{.}
\end{equation}
A simple key observation is that the coboundary operator $\mathrm{d}$ and the Deligne
differential $\mathrm{D}$ commute so that we have a double complex with cochain
groups $C_{\Gamma}^p(\mathcal{D}_{\mathfrak{V}}^k(n))$ in degree $(p,k)$. We denote the total cohomology
of this double complex by $\mathrm{H}^{q}_{\Gamma}(M,\mathcal{D}(n))$. 
In particular, note that we have a natural
group homomorphism
\begin{equation}
\mathrm{eq}:\mathrm{H}^2(M,\mathcal{D}(2)) \to \mathrm{H}_\Gamma^2(M,\mathcal{D}(2)):
\xi \mapsto (\gamma^{*}\xi - \xi,0,0)
\end{equation}
that includes an ordinary Deligne cohomology class into the cohomology of the total
complex we just have defined. 

Now let $\mathcal{G}$ be a hermitian bundle gerbe with connection
and curving with Deligne class $\xi$, and let $(\mathcal{A}_{z},\varphi_{z_1,z_2})$
be a $\Gamma$-equivariant structure on $\mathcal{G}$. For local data $a_{\gamma}\in
C_\Gamma^1(\mathcal{D}_{\mathfrak{V}}^1(2))$ of the isomorphism $\mathcal{A}_\gamma$, i.e.
\begin{equation}
\mathrm{D}a_\gamma = \gamma^{*}\xi-\xi\text{,}
\end{equation}
and local data $b_{\gamma_1,\gamma_2}\in C_\Gamma^2(\mathcal{D}_{\mathfrak{V}}^0(2))$ of
the 2-isomorphisms $\varphi_{\gamma_1,\gamma_2}$, i.e.
\begin{equation}
\mathrm{D}b_{\gamma_1,\gamma_2} = (\mathrm{d}a)_{\gamma_1,\gamma_2}\text{,}
\end{equation}
the commutativity of diagram (\ref{5}) imposes
the condition $(\mathrm{d}b)_{\gamma_1,\gamma_2,\gamma_3}=0$. This means for a bundle
gerbe with Deligne class $\xi$ that the class $\mathrm{eq}(\xi) \in \mathrm{H}^2_\Gamma(M,\mathcal{D}(2))$
is the obstruction class for $\mathcal{G}$ to admit $\Gamma$-equivariant structures. Furthermore,  the cohomology group $\mathrm{H}^1_\Gamma(M,\mathcal{D}(2))$ classifies
the inequivalent choices.

\medskip

In the case of bundle gerbes over a compact, simple, connected and simply-connected
Lie group $G$ we consider the action of a subgroup $Z$ of the center of $G$
by multiplication. In this case the relevant cohomology groups reduce to
the usual group cohomology of the finite group $Z$, so that there is an obstruction
class in $\mathrm{H}^3_{\mathrm{Grp}}(Z,U(1))$, and the possible $Z$-equivariant
structures are classified by $\mathrm{H}^2_{\mathrm{Grp}}(Z,U(1))$. 
For the bundle gerbes we have  constructed above, all obstruction classes
against $Z$-equivariant structures for all subgroups
$Z$ of the center of $G$ can be calculated \cite{gawedzki1,gawedzki2}.

\medskip 
 
Let us now describe how
a choice $(\mathcal{A}_z,\varphi_{z_1,z_2})$
of a $\Gamma$-equivariant structure on a given bundle gerbe\index{bundle
gerbe}
$\mathcal{G}=(\pi,L,\mu,\nabla,C)$ with connection and curving over $M$ defines a quotient bundle gerbe $\overline{\mathcal{G}}=(\overline{\pi},\overline{L},\overline{\mu},\overline{\nabla},\overline{C})$ over $\overline{M}:=M/\Gamma$. Following \cite{gawedzki1}, we
set $\overline{Y}:=Y$ and $\overline{\pi}:=
p \circ
\pi: \overline{Y} \to\overline{M}$, where $p:M \to \overline{M}$ is the projection
to the quotient.
Note that the fibre products
are 
\begin{equation}
\overline{Y}^{[2]}=\bigsqcup_{\gamma\in
\Gamma
}Z^{\gamma}\quad\text{ and }\quad\overline{Y}^{[3]}=\bigsqcup_{\gamma_1,\gamma_2
\in \Gamma} Z^{\gamma_1,\gamma_2}
\end{equation}
for the smooth manifolds 
\begin{equation}
Z^{\gamma} := Y_{\gamma} \times_M Y
\quad\text{ and }\quad
Z^{\gamma_1,\gamma_2}=Y_{\gamma_1\gamma_2}
\times_M Z^{\gamma_1}\text{,}
\end{equation}
 where $Y_{\gamma}:=Y$ as manifolds but with projection
$\gamma^{-1} \circ \pi$ instead of $\pi$. The manifolds $\overline{Y}^{[2]}$
and $\overline{Y}^{[3]}$ have again projections $\overline{\pi}_{i}:=\pi_{i}$
and $\overline{\pi}_{ij}:=\pi_{ij}$ to $\overline{Y}$ and to $\overline{Y}^{[2]}$
respectively.
The curving of the quotient bundle gerbe
will be $\overline{C}:=C \in \Omega^2(\overline{Y})$,
and the line bundle $\overline{L}
\to \overline{Y}^{[2]}$ will be
$L|_{Z^\gamma} := A_{\gamma}$.
Axiom (G1) for the curvature
of $L$ is 
\begin{equation}
\mathrm{curv}(\overline{L})|_{Z^{\gamma}}=\mathrm{curv}(A_{\gamma})=\pi_2^{*}C
- \pi_1^{*}C = \overline{\pi}_2^{*}\overline{C}-\overline{\pi}_1^{*}\overline{C}\text{,}
\end{equation}
and hence satisfied. The isomorphism
$\overline{\mu}$ over $\overline{Y}^{[3]}$
is defined by 
\begin{equation}
\overline{\mu}|_{Z^{\gamma_1,\gamma_2}}:=\varphi_{Z^{\gamma_1,\gamma_2}}
\label{26}\end{equation}
and its associativity by means
of axiom (G2) is nothing
but the condition on
the isomorphisms $\varphi_{\gamma_1,\gamma_2}$
from Definition \ref{def1}. So we have defined a
bundle gerbe $\overline{\mathcal{G}}$
over $\overline{M}$ with connection and curving.

Applying this procedure to the bundle gerbes over the simply connected Lie
groups and their equivariant structures, one obtains examples of bundle
gerbes over all compact simple connected and simply-connected Lie groups.

\section{Structure on Loop Spaces from Bundle Gerbes}

\label{sec1}

We describe a construction of a line bundle over the loop space $LM$ \index{loop
space}of a
manifold $M$ from a given hermitian bundle gerbe  over $M$ with connection
and curving.
This construction is adapted from the one in \cite{brylinski1}. 
 
For preparation, we have to describe the set of isomorphisms between two
fixed hermitian bundle gerbes $\mathcal{G}_1$ and $\mathcal{G}_2$ with connection
and curving. We denote the set of isomorphism classes of morphisms $\mathcal{A}:\mathcal{G}_1
\to \mathcal{G}_2$ by
$\mathrm{Iso}(\mathcal{G}_1,\mathcal{G}_2)$. The following lemma can easily
be shown using the cohomological description.

\begin{lemma}[\normalfont \cite{schreiber1}]
\label{prop1}
The group $\mathrm{Pic}_0^{\nabla}(M)$ of isomorphism classes of flat line bundles
with  connection over $M$ acts freely and transitively on the set $\mathrm{Iso}(\mathcal{G}_1,\mathcal{G}_2)$. \end{lemma}

Now let 
$LM:=C^{\infty}(S^1,M)$ be the free loop space of $M$, equipped with a smooth
manifold structure as described in \cite{brylinski1}. Let  $\mathcal{G}=(\pi,L,\mu,\nabla,C)$
be a hermitian bundle gerbe on $M$ with connection and curving. 
The total space of the line bundle \index{line bundle!over the loop space}we are going to define is, as a set, 
\begin{equation}
\mathcal{L}:= \bigsqcup_{\gamma\in LM} \mathrm{Iso}(\gamma^{*}\mathcal{G},\mathcal{I}_0)\text{.}
\end{equation}
It comes with the evident projection to $LM$, and by Lemma \ref{prop1} every
fibre is a torsor over the group
\begin{equation}
\mathrm{Pic}^{\nabla}_0(S^1)\cong\mathrm{Hom}(\pi_1(S^1),U(1))\cong U(1)\text{.}
\end{equation}
Note that the canonical projection $p: \mathcal{L} \to LM$
admits local sections: for a contractible open
subset $U \subset M$ we have an isomorphism $\mathcal{T}: \mathcal{G}|_U \to \mathcal{I}_{\rho}$ which provides a section $\sigma: LU \to \mathcal{L}: \gamma
\mapsto (\gamma,\gamma^{*}\mathcal{T})$. 

\begin{proposition}
There is a unique differentiable structure on $\mathcal{L}$, such that the
projection $p$ and  the sections $\sigma$ are smooth, and $\mathcal{L}$ becomes
a principal $U(1)$-bundle over $LM$. 
\end{proposition}

\proof
Since every gerbe over $S^1$ is trivializable, none of the fibres $p^{-1}(\gamma)$
is empty. Hence each fibre is  a $U(1)$-torsor. For
the same reason, the
image $\gamma(S^1)$ of any loop $\gamma$
has an open neighbourhood $U \subset M$ (cf. the proof of Proposition 6.2.1 in \cite{brylinski1}),
such that $\mathcal{G}|_U$ admits a trivialization
$\mathcal{T}: \mathcal{G}|_U \to \mathcal{I}_{\rho}$.
The corresponding section $\sigma$ identifies
$p^{-1}(LU)$ with $LU \times U(1)$, and thus
defines a topology and a differentiable structure
on each preimage $p^{-1}(LU)$. Let a topology
on $\mathcal{L}$ be generated by all open subsets of
all the fibres $p^{-1}(LU)$. Now, for two
intersecting subsets $U_1$ and $U_2$ and trivializations
$\mathcal{T}_1$ and $\mathcal{T}_2$ respectively,
let $N$
be the line bundle over $U_1 \cap U_2$ from Lemma \ref{lem1}. The transition map $LU_1
\times U(1) \to LU_2 \times U(1)$ is then
given by $(\gamma,z) \mapsto (\gamma,z \cdot
\mathrm{hol}_{\gamma^{*}N}(S^1)^{-1})$, and hence differentiable with respect to
the loop $\gamma:S^1 \to U_1 \cap U_2$.
\endofproof

Instead of a principal $U(1)$-bundle, we will often and equivalently consider
$\mathcal{L}$ as a hermitian line bundle. The construction just discussed
also applies to the case when the bundle gerbe $\mathcal{G}$ is hermitian and has
a connection. A curving defines a connection $\nabla$ on the line bundle $\mathcal{L}$, whose curvature is
\begin{equation}
\mathrm{curv}(\nabla) = \int_{S^1}\mathrm{ev}^{*}H\text{,}
\end{equation}
where $H$ is the curvature of the gerbe $\mathcal{G}$, and $\mathrm{ev}:LM
\times S^1 \to M$ is the evaluation map. The hermitian line bundle $\mathcal{L}$  with connection is
 functorial in the following sense: 
\begin{proposition}
\label{prop2}
For a hermitian bundle gerbe $\mathcal{G}$ over $M$ with connection and curving, denote
the associated hermitian line bundle with connection  by 
$\mathcal{L}_{\mathcal{G}}$. \begin{itemize}
\item[i)]
Any isomorphism of gerbes $\mathcal{G} \to
\mathcal{G}'$ induces an isomorphism $\mathcal{L}_{\mathcal{G}'}
\to \mathcal{L}_{\mathcal{G}}$ of line bundles.

\item[ii)]
For a smooth map $f:X \to M$ denote the induced
map on loop spaces by $Lf: LX \to LM$. Then, the line bundles  $(Lf)^{*} \mathcal{L}_{\mathcal{G}}$ and $\mathcal{L}_{f^{*}\mathcal{G}}$
are canonically isomorphic. 

\item[iii)]
For the dual gerbe $\mathcal{G}^{*}$ we obtain a canonical isomorphism $\mathcal{L}_{\mathcal{G}}^{*}
\cong \mathcal{L}_{\mathcal{G}^{*}}$.

\end{itemize}
\end{proposition}

\medskip

Let us also have a view on the cohomological counterpart of the construction
of a line bundle over $LM$ from a bundle gerbe over $M$. For this purpose,
one has to extend the usual fibre integration of differential forms,
\begin{equation}
\int_{S_1}: \Omega^{k+1}(X \times S^1) \to \Omega^k(X)\text{,}
\end{equation}
to the Deligne cohomology groups. Here, $X$ can be any smooth and possibly
infinite dimensional manifold. Such extensions have been described in various
ways \cite{gawedzki3, brylinski1,gomi2}. Then, for $X=LM$, the concatenation
of this extension  with the pullback
along the evaluation
map
\begin{equation}
\mathrm{ev}: LM \times S^1 \to M
\end{equation}
gives a group homomorphism
\begin{equation}
\int_{S^1} \circ\, \mathrm{ev}^{*}: \mathrm{H}^{3}(M,\mathcal{D}(3)) \to
\mathrm{H}^2(LM,\mathcal{D}(2))\text{.}
\end{equation}
One can show that the image of the Deligne class $[(g,A,B)]$ of a hermitian
bundle gerbe $\mathcal{G}$ with connection and curving under this group homomorphism
gives exactly the Deligne class of the line bundle $\mathcal{L}$ with the
connection $\nabla$ we have constructed above in a direct geometric way.

\medskip

The construction of the line bundle $\mathcal{L}$ can in particular be applied to the bundle
gerbes over compact Lie groups\index{Lie group!compact} from Section \ref{sec2}. In this case, we obtain a  sequence
\begin{equation}
\label{15}
\alxy{1 \ar[r] & U(1) \ar[r] & \mathcal{L} \ar[r]^-{p} & LG \ar[r] & 1}
\end{equation}
of smooth maps, where $U(1)$ is mapped to the fibre of $\mathcal{L}$ over the loop which
is constantly $1\in G$. Now it is a natural question, whether one can equip the
total space $\mathcal{L}$ with a group structure, such that the sequence
(\ref{15}) is an exact sequence of groups. This would provide a geometric
construction of central extensions\index{central extension} of loop groups. For simply-connected groups
there exist definitions of group structures on $\mathcal{L}$ \cite{brylinski1}, \index{Lie group!simply-connected}while in the
general case the geometric construction of loop group extensions from bundle gerbes  is still an open problem.

\section{Algebraic Structures for Gerbes}

\label{sec5}

There are several additional structures for bundle gerbes, some of which we introduce in this section. We describe the particular  case of those additional
structures on bundle gerbes over compact connected Lie groups.

\subsection{Bundle Gerbe Modules}

Bundle gerbe modules, also known as twisted vector bundles, have been introduced in \cite{bouwknegt1} in order to
realize twisted K-theory geometrically. They are also the appropriate structure
to extend the definition of holonomy to surfaces with boundary \cite{carey2}.

\begin{definition}
\label{def4}
Let $\mathcal{G}$ be a bundle gerbe over $M$. A  \emph{$\mathcal{G}$-module}\index{bundle
gerbe module}
is a 1-morphism $\mathcal{E}:
\mathcal{G} \to \mathcal{I}_{\omega}$ for some 2-form $\omega\in\Omega^2(M)$. The 2-form
$\omega$ is called the curvature of the gerbe module.
\end{definition}

Let us compare this definition with the original definition
of  bundle gerbe modules in \cite{bouwknegt1}. A left
$\mathcal{G}$-module $\mathcal{E}:\mathcal{G} \to \mathcal{I}_{\omega}$ consists of a vector bundle $E$ over $Y$ and of an isomorphism $\epsilon:
L \otimes \pi_2^{*}E \to \pi_1^{*}E$ of vector bundles over $Y^{[2]}$ which satisfies
\begin{equation}
\pi_{13}^{*}\epsilon\circ (\mu \otimes \id)= \pi_{23}^{*}\epsilon \circ \pi_{12}^{*}\epsilon\text{.}
\end{equation}
The similarity with an action of $L$ on $E$ justifies the notion gerbe module.
The curvature of $E$ is restricted by (\ref{23})  to 
\begin{equation}
\frac{1}{n}\mathrm{tr}(\mathrm{curv}(E)) = \pi^{*}\omega - C
\end{equation}
 with $n$ the rank of $E$. 

In terms of local data, a rank-$n$ bundle gerbe module $\mathcal{E}{:}\ \mathcal{G}\to
\mathcal{I}_{\omega}$ is  described by a collection
$(G_{ij},\Pi_i)$ of smooth functions $G_{ij}: U_i \cap U_j \to U(n)$ and $\mathfrak{u}(n)$-valued 1-forms 
$\Pi_i \in \Omega^1(U_i) \otimes \mathfrak{u}(n)$ which relate the local 
data of the bundle gerbes $\mathcal{G}$ and $\mathcal{I}_{\omega}$
in the following way:
\begin{equation}
  \begin{array}{lcll}
  1 &=& g_{ijk}^{} \,\cdot\, G_{ik}^{}\,G_{jk}^{-1}\,G_{ij}^{-1}
  &~\text{ on }~U_i \cap U_j \cap U_k \,,
  \\{}\\[-.7em]
  0 &=& A_{ij}^{} + \Pi_{j}^{} - G_{ij}^{-1}\Pi_{i}^{}\,G_{ij}^{}-
  G_{ij}^{-1}\,\mathrm{d}G_{ij}^{}
  &~\text{ on }~U_i \cap U_j \,,
  \\{}\\[-.9em]
  \omega &=& \,B_i^{}\, + \mbox{$\displaystyle \frac{1}{n}$}\, 
  \mathrm{tr}(\mathrm{d}\Pi_i^{})&~\text{ on }~U_i \,\text{.}\end{array}
\end{equation}
Note that the derivative of the last equality gives
\begin{equation}
\label{19}
\mathrm{d}\omega = \mathrm{d}B_i = \mathrm{curv}(\mathcal{G})\text{.}
\end{equation} 
Also note that if the bundle gerbe $\mathcal{G}$ is itself
trivial, i.e.\ has local data $(1,0,C|_{U_i})$ for a globally defined
2-form $C \in \Omega^2(M)$, then $(G_{ij},\Pi_{i})$ are the
local data of a rank-$n$ vector bundle over $M$ with curvature of trace 
$n\,(\omega{-}C)$. This explains the terminology
\textquotedblleft twisted\textquotedblright\ vector bundle in the non-trivial
case.  

\medskip

According to (\ref{19}), a necessary condition for the existence
of a bundle gerbe module is that the curvature
is an exact form. However, this is not the case in many situations,
for example for the bundle gerbes on compact Lie groups we have constructed
in Section \ref{sec2}, whose curvature is the canonical 3-form $H$. For this reason, one often considers a pair $(Q,\mathcal{E})$ of a submanifold $Q \subset M$ together
with a gerbe module $\mathcal{E}:\mathcal{G}|_{Q} \to \mathcal{I}_{\omega}$
for the restriction of the gerbe to this submanifold. In conformal field
theory, the pair $(Q,\mathcal{G})$ is also called a D-brane.

In particular we can consider this situation for the bundle gerbes over Lie
groups $G$ constructed in Section \ref{sec2}. In this case, the important submanifolds
are conjugacy classes\index{conjugacy class} $Q=\mathcal{C}_\lambda$, and we already know that the
curvature $\mathrm{curv}(\mathcal{G})=H$ becomes exact when restricted to
a conjugacy class, $H|_{\mathcal{C}_{\lambda}}=\mathrm{d}\omega_{\lambda}$.
So the necessary condition (\ref{19}) is satisfied. One can furthermore show
\cite{gawedzki4} that precisely for integrable weights $\lambda$ there exists a $\mathcal{G}|_{\mathcal{C}_{\lambda}}$-module
with curvature $\omega_{\lambda}$.
This is the  appropriate description of \textquotedblleft flux stabilization
of D-branes\textquotedblright\ in string theory \cite{bachas2}.

\subsection{Bundle Gerbe Bimodules}

Bundle gerbe bimodules generalize bundle gerbe modules for one bundle gerbe
$\mathcal{G}$ to a structure for two bundle gerbes.

\begin{definition}[\normalfont \cite{fuchs4}]
Let $\mathcal{G}_1$ 
and $\mathcal{G}_2$ be bundle gerbes over $M$. A \emph{$\mathcal{G}_1$-$\mathcal{G}_2$-bimodule}\index{bundle
gerbe bimodule} 
is a morphism
\begin{equation}
  \mathcal{D} :
  \mathcal{G}_1^{} \to \mathcal{G}_2
  \otimes \mathcal{I}_{\varpi}
\end{equation}
for some 2-form $\varpi \in \Omega^2(M)$. The 2-form $\varpi$  is called the curvature of the bimodule.
\end{definition}
This definition is related to the one of a gerbe module in 
the sense that -- using the appropriate notion of duality for bundle gerbes
 \cite{waldorf1} -- a $\mathcal{G}_1$-$\mathcal{G}_2$-bimodule
is the same as a $(\mathcal{G}_1{\otimes}\mathcal{G}_2^{*})$-module. 

We drop the discussion of the local data of a gerbe bimodule. However, it
is clear that there is -- analogously to equation (\ref{19}) -- a necessary condition on the curvature of the two bundle
gerbes,
\begin{equation}
\label{22}
H_1 = H_2 + \mathrm{d}\varpi\text{.} 
\end{equation}
It is again useful to consider bimodules for restrictions of the bundle gerbes
to a submanifold $Q \subset M$. In particular, if $M$ is the direct product
of two manifolds $M_1$ and $M_2$, each equipped with a bundle gerbe $\mathcal{G}_1$
and $\mathcal{G}_2$ respectively,  one considers   $p_1^{*}\mathcal{G}_1|_Q$-$p_2^{*}\mathcal{G}_2|_Q$-bimodules,
where $p_i:M \to M_i$ is the projection on the $i$th factor.
This is the setup to define holonomies around surfaces with defect lines
\cite{fuchs4}.

Examples for such bimodules are again provided by compact Lie groups and
the basic bundle gerbes thereon \cite{fuchs4}. The study of such examples
leads to the following relevant submanifolds
$Q$
of $G \times G$, so-called biconjugacy classes \index{biconjugacy class}
\begin{equation}
  \mathcal{B}_{h_1,h_2} := \left\lbrace (x_1h_1x_2^{-1},x_1h_2x_2^{-1})\in G \times G \;|\;
  x_1, x_2 \in G \right\rbrace 
\end{equation}
for any pair $(h_1,h_2)\in G \times G$.

Biconjugacy classes inherit from the diagonal left and diagonal right actions
of $G$ on $G\times G$ two commuting actions of $G$.  One observes that the smooth map
\begin{equation}
  \tilde \mu :  G \times G \to G: (g_1,g_2) \mapsto g_1g_2^{-1}
\end{equation}
intertwines the diagonal left and diagonal right action of
$G$ on $G\times G$ and the adjoint and trivial actions of $G$ on itself,
respectively. It now follows that a biconjugacy class in $G \times G$ is the preimage of 
a conjugacy class in $G$ under the  projection $\tilde \mu$:
\begin{equation}
  \mathcal{B}_{h_1,h_2}^{} = \tilde\mu^{-1}(\mathcal{C}_{h_1^{}h_2^{-1}}) =
  \big\lbrace (g_1,g_2) \in G \times G \;|\; g_1^{}g_2^{-1}
  \in \mathcal{C}_{h_1^{}h_2^{-1}} \big\rbrace \,\text{.}
\end{equation}
We introduce the two-form 
\begin{equation}
  \varpi_{h_1^{},h_2^{}}:= \tilde\mu^* \omega_{h_1^{}h_2^{-1}} -
  \frac{k}{2}\,\langle p_1^*\theta\wedge p_2^*\theta\rangle
\label{3.15}
\end{equation}
on $\mathcal{B}_{h_1,h_2}$, where 
both summands are restricted to the submanifold $\mathcal{B}_{h_1,h_2}$
of $G\times G$. From the intertwining properties of $\tilde\mu$ and the bi-invariance
of $\omega$ it follows
that the two-form $\varpi$ is also bi-invariant. One can  show that
it satisfies 
\begin{equation}
  p_1^* H = p_2^* H +\mathrm{d}\varpi_{h_1,h_2}
\end{equation}
on a biconjugacy class
$\mathcal{B}_{h_1,h_2}$, which is the necessary condition (\ref{22}) \cite{fuchs4}.

\subsection{Jandl Structures}

Another structure we want to introduce is a Jandl structure on a bundle gerbe
$\mathcal{G}$. Jandl structures extend the definition of holonomy we gave
in Section \ref{sec7} to unoriented, and in particular  unorientable surfaces.

\begin{definition}[\normalfont \cite{schreiber1}]
A \index{Jandl structure}\emph{Jandl structure} $\mathcal{J}$ on a bundle gerbe $\mathcal{G}$ over $M$ is an involution $k$ of $M$ together with an isomorphism 
\begin{equation}
\mathcal{A}:k^{*}\mathcal{G} \to \mathcal{G}^{*}
\end{equation}
and a 2-morphism
\begin{equation}
\varphi: k^{*}\mathcal{A} \Rightarrow \mathcal{A}^{*}
\end{equation}
which satisfies the equivariance condition
\begin{equation}
k^{*}\varphi = \varphi^{*}\text{.}
\end{equation}
\end{definition}

To give  an impression of the details of a Jandl structure, recall
that an isomorphism such as $\mathcal{A}=(A,\alpha)$ consists of a line bundle $A$ over
the space $Z$ which is build up from the two surjective submersions of the bundle gerbes
$k^{*}\mathcal{G}$ and $\mathcal{G}^{*}$. In this particular situation, there is a canonical lift $\tilde k$ of the involution $k$ into the space $Z$, and it is in fact
easy to work out that the 2-morphism $\varphi$
 defines a $\tilde k$-equivariant structure
on the line bundle $A$, which is compatible with the isomorphism $\alpha$. Summarizing, a Jandl structure $\mathcal{J}$ on $\mathcal{G}$ is
an isomorphism
\begin{equation}
\mathcal{A}:k^{*}\mathcal{G} \to \mathcal{G}^{*}
\end{equation} 
whose line bundle $A$ is equivariant with respect to the involution
$\tilde k$ on $Z$ \cite{schreiber1}. 

\medskip 

Recall that we introduced equivariant structures \index{equivariant structure!on
a bundle gerbe} on bundle gerbes in Section
\ref{sec2} in order to produce bundle gerbes over quotients of a manifold
$M$ by a discrete group $Z$. One can combine equivariant structures and
Jandl structures to $Z$-equivariant Jandl structures, leading to the mathematically
appropriate description of so-called orientifolds in string theory \cite{gawedzki6}.  The idea behind
this combination is, that a bundle gerbe $\mathcal{G}$ over $M$ with $\Gamma$-equivariant Jandl
structure defines a bundle gerbe $\overline{\mathcal{G}}$ with Jandl structure  over the quotient
$M/Z$. For a cohomological description, one modifies the action of $Z$ on
the Deligne cohomology group to an action of the semidirect product $\Gamma := \mathbb{Z}_2 \ltimes Z$ \cite{gawedzki6}. 

For bundle gerbes over compact Lie groups, where $Z$ is a subgroup of the
center of $G$, the relevant involutions are given by
\begin{equation}
k_z: G \to G : g \mapsto (zg)^{-1}
\end{equation}
for any $z$ in the center. Again, using the basic bundle gerbes constructed
in Section \ref{sec2}, one can classify all equivariant Jandl structures
over all these bundle gerbes \cite{gawedzki6}.

\section{Applications  to Conformal Field Theory}

\label{sec6}

Let us explain the relation between  conformal field theory and Lie
theory, which  arises in the study of non-linear sigma models on a Lie group
$G$. Such a model can be defined
by  amplitudes $\mathcal{A}(\phi)$ for some
path integral, where  $\phi$ is a map from a closed complex curve $\Sigma$
-- the world sheet -- into the target space $G$ of the model. In \cite{witten1}, Witten gives the following definition for $G=SU(2)$.
$\Sigma$ is the boundary of a three dimensional manifold
$B$, and because the homotopy groups $\pi_i(SU(2))$ vanish
for $i=1,2$, every map
$\phi:\Sigma \to M$ can be extended into the interior $B$ to a map $\Phi:B \to G$. Witten showed that -- due to the integrality
of  the canonical 3-form $H$ --
\begin{equation}
\label{17}
\mathcal{A}(\phi):=\exp \left (S_{\mathrm{kin}}(\phi) + \int_B \Phi^{*}H \right ) 
\end{equation}
neither depends on the choice of $B$ nor on the choice of the extension
$\Phi$, so that one obtains a well-defined amplitude. Here $S_{\mathrm{kin}}(\phi)$
is a kinetic term, and with a certain relative normalization of the two terms
in (\ref{17}) this
model is called the  \index{Wess-Zumino-Witten model}Wess-Zumino-Witten
model on $G$ at level $k$. For  non-simply-connected Lie groups, the extension $\Phi$ of $\phi$ to $B$ does not
exist in general. In these cases, the second summand of (\ref{17}) has to be generalized.

\begin{proposition}
\label{prop3}
Let $\mathcal{G}$ be a hermitian bundle gerbe over
$G$ with connection and curving of  curvature $H$. For a three-dimensional oriented manifold $B$ with boundary and a map $\Phi:B
\to G$, we have 
\begin{equation}
\mathrm{hol}_{\mathcal{G}}(\Phi|_{\partial B}) = \exp \left ( \int_B \Phi^{*}H
\right )\text{.}
\end{equation} 
\end{proposition}

\proof 
Remember that for any trivialization $\mathcal{T}:
\Phi^{*}\mathcal{G}|_{\partial B} \to \mathcal{I}_{\rho}$
we have $\Phi^{*}H|_{\partial B}=\mathrm{d}\rho$.
The rest follows by Stokes' Theorem.  
\endofproof

This way we reproduce the amplitude
of the coupling term 
of the Wess-Zumino-Witten model
\index{Wess-Zumino-Witten model}by
\begin{equation}
\mathcal{A}(\phi) = \exp \left
( S_{\mathrm{kin}}(\phi) \right
) \cdot \mathrm{hol}_{\mathcal{G}}(\phi)\text{.}
\end{equation}
Notice that using bundle gerbes we did not impose any
condition on the topology of the
target space $G$. For compact connected and simply-connected Lie
groups, it reproduces Witten's original definition. However, for general target spaces
there may be bundle gerbes with
same curvature, which are not isomorphic.
This occurs 
for instance  for the Lie group $\mathrm{Spin}(4n)/(\mathbb{Z}_2
\times \mathbb{Z}_2)$. Here, the theory of bundle gerbes over Lie groups
has brought new insights into the Lagrangian description of Wess-Zumino-Witten
models.

\medskip

In conformal field theory  many applications require to consider surfaces with
boundary. For those, we are not
able to apply Definition
\ref{def2} of the holonomy of a bundle
gerbe $\mathcal{G}$: for a change of the chosen trivialization, a boundary term emerges which has to be
compensated to achieve a  holonomy
 independent of the choice of the trivialization. The compensating term is
 provided by the choice of a symmetric D-brane $(\mathcal{C}_{\lambda},\mathcal{E}_{\lambda})$:
 \index{D-brane} a conjugacy class $\mathcal{C}_{\lambda}$ for an integrable weight $\lambda$,
together with a $\mathcal{G}|_{\mathcal{C}_{\lambda}}$-module $\mathcal{E}_{\lambda}$
of curvature $\omega_{\lambda}$  \cite{gawedzki4}.

Another class of conformal field theories involves  unoriented world sheets $\Sigma$.
Again, Definition \ref{def2} has to be generalized, since it involves the
integral of a differential form over $\Sigma$.  It has been shown in \cite{schreiber1}
that the choice of a Jandl structure on the bundle gerbe $\mathcal{G}$ makes the holonomy again well-defined.
The precise classification of Jandl structures leads to a complete classification of unoriented Wess-Zumino-Witten models \cite{schreiber1,gawedzki6}. 

\medskip

Let us also indicate the relevance of the line bundle $\mathcal{L}$ over
the loop group $LG$ we have constructed in Section \ref{sec1}. The Hilbert
space of holomorphic sections in $\mathcal{L}$ (completed with respect to
its hermitian metric) serves as the space of states for the quantized theory
\cite{gawedzki1}. The choice of additional structures like bundle gerbe modules
or Jandl structures, has implications on this space. 
For example,
a Jandl structure on the bundle gerbe $\mathcal{G}$ implies by Proposition
\ref{prop2} an isomorphism $\varphi: Lk^{*}\mathcal{L} \to \mathcal{L}^{*}$
which satisfies $Lk^{*}\varphi=\varphi$.

\section{Open Questions}

We conclude this contribution with the discussion of  some lines
for further research. One obvious direction is to extend the results
explained here to gauged Wess-Zumino-Witten models, so-called
coset theories. Those models having fixed points under the action of the 
group implementing field identifications for which the so-called untwisted stabilizer 
is strictly smaller than the stabilizer \cite{fuchs6} should be
particularly interesting: in this case, simple gerbe modules
and bimodules of rank strictly bigger than one appear naturally.
A precise understanding of such theories requires the notion of an
$H$-equivariant gerbe (bi-)module on the ambient group $G$.

Another subtle issue is the generalization of our results to non-compact 
Lie groups; this is partially due to the fact that much less is known
about these theories in algebraic approaches.

The following two points seem to be conceptually appealing questions:
our initial motivating question in this paper was about central extensions 
of Lie groups. While hermitian bundle gerbes naturally account for a line 
bundle $\mathcal{L}$ on loop space, more specific structure on the gerbe is 
needed to obtain a group structure on the line bundle $\mathcal{L}$. 
Similarly, one wishes to find structure on a gerbe module $E$ that endows 
the associated bundles $\mathcal{E}$ over loop- and interval spaces with a 
natural $\mathcal{L}$-module structure.

Finally, we point out that gerbe bimodules have
a natural operation of fusion which is very much in spirit of
a convolution of correspondences. Imposing additional properties
on gerbe bimodules over a compact Lie group (that in physical applications 
ensure the existence of enough conserved quantities) one should be able to 
single out interesting subcategories of gerbe bimodules that are semi-simple 
tensor categories. If $G$ is simply connected, their fusion ring can be 
expected to be the
corresponding Verlinde algebra; for non-simply connected groups,
we expect interesting cousins of the Verlinde algebra.

\newcommand{\etalchar}[1]{$^{#1}$}

\end{document}